\documentclass[final]{amsart}

\usepackage{xspace}
\usepackage[matrix,arrow,curve]{xy}
\usepackage[backref,bookmarks,breaklinks,colorlinks,unicode]{hyperref}

\newcommand{\Pre}[1]{\widehat{#1}}
\newcommand{\CoPre}[1]{\check{#1}}
\newcommand{\Pro}[2][]{\operatorname{Pro}_{#1}({#2})}
\newcommand{\Ind}[2][]{\operatorname{Ind}_{#1}({#2})}
\newcommand{\y}{\mathbf{y}}
\newcommand{\yOp}{\mathbf{\check{y}}}
\newcommand{\CC}{\mathcal{C}}
\newcommand{\DD}{\mathcal{D}}
\newcommand{\MM}{\mathcal{M}}
\newcommand{\Lim}[2][]{\underset{\stackrel{\longleftarrow}{#1}}{\lim}\,#2}
\newcommand{\CoLim}[2][]{\underset{\stackrel{\longrightarrow}{#1}}{\lim}\,#2}
\newcommand{\Op}[1]{{{#1}^\circ}}
\newcommand{\Hom}[3][]{\operatorname{Hom}_{#1}(#2,#3)}
\newcommand{\Def}[1]{\emph{#1}\xspace}
\newcommand{\ra}{\rightarrow}
\newcommand{\x}{\times}
\newcommand{\p}{\mathfrak{p}}
\newcommand{\inv}[1]{{#1}^{-1}}
\newcommand{\w}{\mbox{${\omega}$}}
\newenvironment{Thanks}{\subsection*{Acknowledgement}}{}

\newcommand{\newthm}[2]{\newtheorem{#1}[lemma]{#2}}
\theoremstyle{plain}
\newtheorem{lemma}{Lemma}
\newthm{Ass}{Assumption}
\newthm{convention}{Convention}
\newthm{criterion}{Criterion}
\newthm{fact}{Fact}
\newthm{cor}{Corollary}
\newthm{claim}{Claim}
\newthm{thm}{Theorem}
\newthm{prop}{Proposition}

\theoremstyle{definition}
\newthm{definition}{Definition}
\newthm{example}{Example}

\theoremstyle{remark}
\newthm{remark}{Remark}
\newthm{notation}{Notation}

\newcommand{\DefAlias}[2]{\expandafter\xdef\csname #1\endcsname{#2}}
\newcommand{\CiteAlias}[2]{\DefAlias{CITE#1}{#2}}
\newcommand{\Cite}[1]{\cite{\csname CITE#1\endcsname}}

\CiteAlias{ACFA}{MR1652269}
\CiteAlias{LAG}{MR1102012}
\CiteAlias{DG}{MR1972449}
\CiteAlias{groupoid}{groupoid}
\CiteAlias{bounds}{MR739626}
\CiteAlias{stability}{MR1416106}
\CiteAlias{SGA}{MR0354652}
\CiteAlias{singer}{MR1480919}
\CiteAlias{milne}{MR559531}
\CiteAlias{topoi}{MR1300636}
\CiteAlias{sacks}{MR0398817}
\CiteAlias{EI}{MR1083551}
\CiteAlias{HrWeil}{MR1827833}
\CiteAlias{SingerMod}{SingerMod}
\CiteAlias{PillayZoe}{MR1650667}
\CiteAlias{Bouscaren}{MR984628}
\CiteAlias{eisenbud}{MR1322960}

\title{Ind- and Pro- definable sets}
\author{Moshe Kamensky}

\address{Department of Mathematics\\
         The Hebrew University\\
         Jerusalem, Israel}
\curraddr{Department of Maths\\
         University of East Anglia\\
         Norwich, NR4 7TJ, England}

\newcommand{\Email}{\url{mailto:m.kamensky@uea.ac.uk}}
\email{\Email}


\begin{document}

\bibliographystyle{amsplain}

\begin{abstract}
  We describe the ind- and pro- categories of the category of definable sets, 
  in some first order theory, in terms of points in a sufficiently saturated 
  model.
\end{abstract}

\maketitle

\section{Introduction}
Given the direct limit $Y$ of some system $Y_i$ in a given category, the 
morphisms from $Y$ to another object $X$ are described, by definition, as 
certain collections of morphisms from each $Y_i$ to $X$. In contrast, there 
is, in general, no simple description of morphisms in the other direction, 
from $X$ to $Y$. However, if the category in question is, for example, a 
category of topological spaces, and $X$ is compact, then any morphism from 
$X$ to $Y$ will factor via some $Y_i$.

The category $\Ind{\CC}$ of ind-objects of a category $\CC$ is a category 
containing the original category $\CC$, in which any filtering system has a 
limit, and the objects of the original category are ``compact'' in the above 
sense. This construction, which appears in \Cite{SGA}, can be applied to any 
category, and is described below. The dual construction, of the category of 
pro-objects, is described as well.

In the context of first order logic, and definable sets, there is a natural 
notion of compactness, and given a system of definable sets, one may compute 
limits of their points in a given model. The purpose of this note is to 
describe how the categorical notions of ind- and pro- objects apply to 
definable sets, and in particular to describe the categories of ind- and pro- 
definable sets in terms of points in a model. The main results are 
proposition~\ref{prp:points}, which explains how to compute the $M$ points of 
$P$, where $M$ is any model and $P$ is an ind-definable (or pro-definable) 
set, and proposition~\ref{prp:morphisms}, which describes morphisms in terms 
of such points. The final statement of the results is in 
corollary~\ref{cor:proind}.

\begin{Thanks}
  This work is part of my PhD research, performed in the Hebrew university 
  under the supervision of Ehud Hrushovski. I would like to thank him for his 
  guidance.
\end{Thanks}

\section{Categorical notions}

We begin by recalling some general notions from category theory. The 
reference to all this is \Cite{SGA}. Let $\CC$ be a category (which we assume 
to be small), $\Pre{\CC}$ the category of presheaves on $\CC$ (i.e., 
contra-variant functors from $\CC$ to the category of sets), and 
$\y:\CC\rightarrow\Pre{\CC}$ the Yoneda embedding, given by 
$\y(X)(Z)=\Hom{Z}{X}$.

A \Def{filtering category} is a small category $I$ such that:

\begin{itemize}
    
  \item For any two objects $i,j$ of $I$, there are morphisms $i\ra{}k$ and 
    $j\ra{}k$ for some object $k$.
    
  \item For any two morphism $t_1,t_2:i\ra{}j$ there is a morphism 
    $s:j\ra{}k$ with $s\circ{}t_1=s\circ{}t_2$.
    
\end{itemize}

A \Def{filtering system} in $\CC$ is a functor from a filtering category to 
$\CC$. Such a system will be denoted $(X_i)$, where $X_i$ is the object of 
$\CC$ associated with $i$. We now define $\Ind{(X_i)}$, the \Def{ind-object} 
of $\CC$ associated with the system $(X_i)$, to be $\CoLim{\y(X_i)}$ (an 
object of $\Pre{\CC}$.) Recall that direct limits in $\Pre{\CC}$ can be  
computed ``pointwise''. Thus, we have for any object $Y$ of $\CC$,
\begin{equation*}
  \Hom{Y}{\Ind{(X_i)}}=\Ind{(X_i)}(Y)=\CoLim{\Hom{Y}{X_i}}
\end{equation*}

The category $\Ind{\CC}$ is defined to be the full subcategory of $\Pre{\CC}$ 
of presheaves isomorphic to $\Ind{(X_i)}$ for some filtering system $(X_i)$.

Any directed partially ordered set can be viewed as a filtered category, and 
conceptually a filtering system can be thought of as a partially ordered one.  
In fact, it can be shown that any filtering system is isomorphic to a 
partially ordered one. However, in some cases (such as the proof of 
proposition~\ref{prp:points} below), the natural index category has the more 
general form.

The category of pro-objects $\Pro{\CC}$ is defined by dualising: it is 
defined to be $\Op{\Ind{\Op{\CC}}}$, where $\Op{\CC}$ denotes the opposite 
category to $\CC$. We describe it explicitly in terms of $\CC$ itself: let 
$\CoPre{\CC}=\Pre{\Op{\CC}}$ be the category of co-variant functors from 
$\CC$ to sets, $\yOp:\CC\ra\CoPre{\CC}$ the (contra-variant) Yoneda 
embedding. Given a co-filtering system $(X_i)$ in $\CC$ (i.e., a 
contra-variant functor from a filtering category to $\CC$), the associated 
pro-object is defined to be the functor $\Pro{(X_i)}=\CoLim{\yOp(X_i)}$. For 
any object $Y$ of $\CC$ we get
\begin{equation*}
  \begin{split}
\Hom[{\Pro{\CC}}]{\Pro{(X_i)}}{Y}=\Hom[{\CoPre{\CC}}]{Y}{\Pro{(X_i}})=\\
=\Pro{(X_i)}(Y)=\CoLim{\Hom[{\CC}]{X_i}{Y}}
\end{split}
\end{equation*}

More generally, we have the following formulas for the morphism sets in the 
$Pro$ and $Ind$ categories:

\begin{subequations}\label{eqn:homind}
\begin{gather}
  \Hom{\Ind{(X_i)}}{\Ind{(Y_j)}}=\Lim[i]{\CoLim[j]\Hom{X_i}{Y_j}}\\
  \Hom{\Pro{(X_i)}}{\Pro{(Y_j)}}=\Lim[j]{\CoLim[i]\Hom{X_i}{Y_j}}
\end{gather}
\end{subequations}

It follows that any presheaf $P$ on $\CC$ extends canonically to $\Pro{\CC}$ 
by setting $P(\Pro{(X_i)})=\CoLim{P(X_i)}$: a map of pro-objects
\begin{equation*}
  f:\Pro{(X_i)}\ra\Pro{(Y_j)}
\end{equation*}
is represented by a sequence of maps $f_j:X_{i_j}\ra{}Y_j$, hence we get maps 
$P(f_j):P(Y_j)\ra{}P(X_{i_j})$ that represent a map from $\CoLim{P(Y_j)}$ to 
$\CoLim{P(X_i)}$. Likewise, any functor from $\CC$ to sets can be extended to 
a functor on $\Ind{\CC}$.

Given an object $X$ of $\CC$, the category $\CC/X$ is defined to have 
$\CC$-morphisms $Y\ra{}X$ as objects, and $\CC$-morphisms over $X$ as 
morphisms.  Then $\Ind{\CC/X}=\Ind{\CC}/X$ and $\Pro{\CC/X}=\Pro{\CC}/X$. The 
first assertion follows by definition (and is true when $X$ is replaced by 
any presheaf), while the second uses the fact that the systems are filtered.

We are going to use the following lemma, which describes a sufficient 
condition for a morphism with a section to be an isomorphism:
\begin{lemma} \label{lma:indiso}
  Let $f:\Ind{(X_i)}\ra{}Y$, $g:Y\ra{}X_0$ be two morphisms, such that 
  $f_0\circ{}g$ is the identity on $Y$. Assume that for any $i$, there is a 
  morphism $t_i:X_i\ra{}X_j$ in the system, such that for any two morphisms 
  $h_1,h_2:V\ra{}X_i$, if $f_i\circ{}h_1=f_i\circ{}h_2$, then 
  $t_i\circ{}h_1=t_i\circ{}h_2$ (this is the formal analogue of saying that 
  $f_j$ is injective on the image of $t_i$.)

  Then $f$ is an isomorphism with inverse $g$.
\end{lemma}
\begin{proof}
  First note that for any filtering system $(X_i)$ and an object $X$ in the 
  system, the (full) subsystem consisting of all objects that have a system 
  morphism from $X$ is isomorphic (in the $Ind$ category) to the original 
  one. Thus we may assume that there is a system morphism from $X_0$ to any 
  other object in the system.

  To show that $g$ is the inverse of $f$, we need to show that $g\circ{}f$ is 
  the identity on $\Ind{(X_i)}$ (the other composition is the identity by 
  assumption.) This amounts to showing that for any $i$, $g\circ{}f_i$ is 
  identified with some morphism in the system $(X_i)$. In other words, we 
  need to show that there are morphisms $t:X_i\ra{}X_k$, $s:X_0\ra{}X_k$ such 
  that $s\circ{}g\circ{}f_i=t$ (In fact, for any object $Z$,
  \begin{equation*}
    \begin{split}
    (g\circ f)_Z(\Ind{(X_i)}(Z)) = (g\circ f)_Z(\CoLim[i]\Hom{Z}{X_i})=\\
    =\CoLim[i]\{g\circ f_i\circ u|u\in\Hom{Z}{X_i}\}
    \end{split}
  \end{equation*}
  If the above condition holds, the map taking $u:Z\ra{}X_i$ to 
  $g\circ{}f_i\circ{}u$ is an isomorphism of the limit sets, since 
  $s\circ{}g\circ{}f_i\circ{}u=t\circ{}u$.)
  
  The situation is this:
  \begin{equation*}
    \vcenter
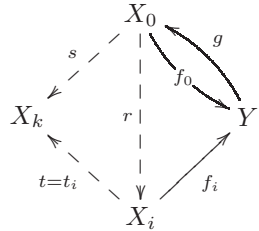
{\xymatrix{
         & X_0\ar@{-->}[dd]_{r}\ar@/_/[dr]|{f_0}\ar@{-->}[ld]_s & \\
     X_k &    & Y\ar@/_/[ul]_{g} \\
         & X_i\ar[ur]_{f_i}\ar@{-->}[ul]^{t=t_i}   &
       }}
  \end{equation*}
  We should find $X_k$, $t$ and $s$, such that the external square commutes. 
  We take $t=t_i$, as promised by the assumption. By the reduction above, 
  there is some morphism $r$ from $X_0$ to $X_i$. We set $s=t\circ{}r$. Thus 
  we should prove that $t\circ{}r\circ{}g\circ{}f_i=t$. By the property of 
  $t$, it is enough to show that $f_i\circ{}r\circ{}g\circ{}f_i=f_i$. But 
  this is true since $f_i\circ{}r\circ{}g=f_0\circ{}g=1_Y$.
\end{proof}

\begin{remark} \label{rmk:indiso}
  In the case that $\CC$ has finite inverse limits, we may replace the 
  arbitrary $V$ by $X_i\x_Y{}X_i$ (and the $h_i$ by the projections.) Thus, 
  in this case we get the following simpler condition:

  Let $\CC$ be a category with finite inverse limits. Let 
  $f:\Ind{(X_i)}\ra{}Y$, $g:Y\ra{}X_0$ be morphisms, such that $f_0\circ{}g$ 
  is the identity on $Y$. Assume that for any $i$, there is a morphism 
  $t_i:X_i\ra{}X_j$ in the system, such that the map 
  $X_i\x_{X_j}X_i\ra{}X_i\x_Y{}X_i$ is an isomorphism.

  Then $f$ is an isomorphism with inverse $g$.
\end{remark}

\begin{remark} \label{rmk:proiso}
  For convenience, we rephrase the above statement in terms of $Pro$ objects:

  Let $\CC$ be a category with finite direct limits. Let $f:Y\ra\Pro{(X_i)}$, 
  $g:X_0\ra{}Y$ be morphisms, such that $g\circ{}f_0$ is the identity on $Y$. 
  Assume that for any $i$, there is a morphism $t_i:X_j\ra{}X_i$ in the 
  system, such that the map $X_i\amalg_{Y}X_i\ra{}X_i\amalg_{X_j}{}X_i$ is an 
  isomorphism.

  Then $f$ is an isomorphism with inverse $g$.
\end{remark}

\section{The case of definable sets}

We now consider the model theoretic setting. The basic terminology is 
explained, for example, in \Cite{sacks}. Let $T$ be a first order theory, 
$\MM$ the opposite category to the category of models of $T$ and elementary 
maps, and $\DD$ the category of definable sets and definable functions 
between them (the word ``definable'' will mean definable over $0$.) The 
relationship between them is described by the faithful functor 
$p:\DD\ra\Pre{\MM}$, given by $p(X)(M)=X(M)$. We first show that this functor 
has a natural extension to the whole category $\Pre{\DD}$ of presheaves on 
$\DD$.

\begin{prop}\label{prp:points}
  There is a fully faithful functor $d:\MM\rightarrow{}\Pro{\DD}$ such that 
  for any definable set $X$ and model $M$, $\Hom{d(M)}{X}=X(M)$. In 
  particular, for any presheaf or functor $F$ on $\DD$, $F(M)$ is well 
  defined.
\end{prop}

Before giving the proof, we roughly explain the idea. A basic property of any 
definable set $X$ is that if $a\in{}X(M)\subseteq{}M^n$, then the whole type 
of $a$ (over $0$) is contained in $X$, and we would like this property to 
hold for an arbitrary presheaf. Since a type is just an example of a 
pro-definable set, $X(M)$ can be written as
\begin{equation*}
  \begin{split}
  X(M)&=\coprod_{a\in M^n}\Hom{tp(a)}{X}=\coprod_{a\in M^n} 
  \CoLim[\substack{\text{$Y$ with}\\
                     a\in Y(M)}]\Hom{Y}{X}=\\
  &=\CoLim[(Y,a\in{}Y(M))]\Hom{Y}{X}
  \end{split}
\end{equation*}
where $\Hom{X}{Y}$ here is taken in the sense of inclusions (so $\Hom{X}{Y}$ 
contains one element if $X\subseteq{}Y$, and is empty otherwise.) When we 
wish to describe this observation in terms of the pro-definable category, we 
run into several problems: first, we obtain distinct systems for distinct 
values of $n$. Second, these systems are not co-filtering.  Finally, it is 
not clear how to distinguish inclusions inside the category.  Fortunately, 
all of these problems are solved by replacing inclusions by arbitrary 
definable maps, as we do in the proof, below.

\begin{proof}[Proof of proposition~\ref{prp:points}]
  Given a model $M$, let $(X_{(X,a)})$ be the system where $a\in{}M$, $X$ is 
  a definable set with $a\in{}X(M)$, and $X_{(X,a)}=X$ (since we no longer 
  distinguish inclusions, we also don't distinguish between elements and 
  tuples.) The morphisms from $X_{(X,a)}$ to $Y_{(Y,b)}$ are definable maps 
  $f:X\rightarrow{}Y$ with $f(a)=b$. This system is cofiltering since all 
  finite inverse limits exist in $\DD$. We abbreviate $X_{(X,a)}$ as $X_a$ 
  and set $d(M)=\Pro{(X_a)}$.  We first show that for any definable set $Y$, 
  we have a canonical bijection $\Hom{d(M)}{Y}\rightarrow{}Y(M)$. Indeed, by 
  definition
  \begin{equation*}
  \Hom{d(M)}{Y}=\CoLim{\Hom{X_a}{Y}}
  \end{equation*}

  So to give a map from $\Hom{d(M)}{Y}$ to $Y(M)$ is the same as to give a 
  matching collection of maps from each $\Hom{X_a}{Y}$ to $Y(M)$. For each 
  $f\in{}\Hom{X_a}{Y}$ we assign $f(a)$. To show that this map is a bijection, 
  we note that the map in the other direction is given by assigning to each 
  $a\in{}Y(M)$ the identity map on $Y=Y_a$. This is, in fact, the inverse, 
  since any definable map $f:X_a\ra{}Y$ is identified with the identity map 
  when ``restricted'' to the graph of $f$. More verbosely, let $f:X_a\ra{}Y$ 
  represent an element in $\Hom{d(M)}{Y}$. Applying the composition of the 
  two maps, we get the identity map on $Y_{f(a)}$. If $\Gamma$ is the graph 
  of $f$, the two projections give maps in the system 
  $\Gamma_{(a,f(a))}\ra{}X_a$ and $\Gamma_{(a,f(a)}\ra{}Y_{f(a)}$ that 
  identify $f$ and the identity on $Y_{f(a)}$.

  To define $d$ on morphisms, we first note that, by what was just shown, 
  given two models $M$ and $N$,
  \begin{equation*}
  \Hom{d(M)}{d(N})=\Lim{\Hom{d(M)}{X_a}}=\Lim{X_a(M)}
  \end{equation*}
  (where the limit is taken over pairs $(X,a)$ with $a\in{}N$.) Thus, to 
  define the map $d:\Hom[{\MM}]{M}{N}\rightarrow{}\Hom{d(M)}{d(N})$ we need to 
  assign, to each elementary map $f:N\rightarrow{}M$ a compatible system of 
  points in the $X_a(M)$. We do it by taking the point $f(a)$. In the other 
  direction, given a matching collection of points, we construct a map from 
  $N$ to $M$ by assigning to a point $a\in{}N$ the point specified for 
  $U_a(M)$ (where $U$ is the universe, $x=x$.) To show that this map is 
  elementary, we note that if, for some definable set $X\subseteq{}U^n$, we 
  have $\bar{a}\in{}X(N)$, $f(\bar{a})$ is the specified point in 
  $X_{\bar{a}}(M)$: $f(a_1,\dots,a_n)=(f(a_1),\dots,f(a_n))$, and each 
  $f(a_i)$ is the specified point for $U_{a_i}$. The projections from 
  ${U^n}_{\bar{a}}$ to the $U_{a_i}$ now show that $f(\bar{a})$ is the 
  specified point for ${U^n}_{\bar{a}}$, hence, because of the inclusion 
  $X\subseteq{}U^n$, for $X_{\bar{a}}(M)$.

  This concludes the construction of the embedding. The last remark follows 
  directly from the remarks above. Explicitly, for $P$ a presheaf on $\DD$, 
  and $F$ a functor from $\DD$ to sets, we have for any model $M$:
  
  \begin{subequations}\label{eqn:points}
  \begin{gather}
    F(M)=\Hom[{\CoPre{\DD}}]{F}{d(M})=\Hom[{\CoPre{\DD}}]{F}{\CoLim{\yOp(X_a}})\\
    P(M)=\CoLim P(X_a)
  \end{gather}
  \end{subequations}
\end{proof}

\begin{remark}
  Instead of viewing definable sets as functors on the category of models, we 
  may, conversely, view a model as a functor on the definable sets. From this 
  point of view, the construction of $d(M)$ (for a general functor $M$) is 
  mentioned in \Cite{topoi} as the \emph{Grothendieck construction}.  
  Unfortunately, I do not know the purpose of this construction in general.
\end{remark}

We are interested in two special cases of the formulas \eqref{eqn:points}: 
let $P=\Ind{(Z_i)}$, $F=\Pro{(Y_i)}$. In this case we obtain:
\begin{gather*}
  \begin{split}
  \Ind{(Z_i)}(M)&=\CoLim{\Ind{(Z_i)}(X_a)}=
  \CoLim[(X,a)]{\CoLim[i]{\Hom{X_a}{Z_i}}}=\\
  &=\CoLim[i]{\CoLim[(X,a)]{\Hom{X_a}{Z_i}}}=
  \CoLim{Z_i(M)}
\end{split}
\\
\begin{split}
  \Pro{(Y_i)}(M)&=\Hom[{\CoPre{\DD}}]{\Pro{(Y_i)}}{d(M})=\\
  &=\Lim[i]{\Hom[{\CoPre{\DD}}]{Y_i}{d(M})}=\Lim[i]{Y_i(M)}
\end{split}
\end{gather*}

Thus, to compute the points of a pro-definable set in a model $M$, we need to 
choose a presentation of it as system, and compute the inverse limit of the 
associated system of sets (and similarly for ind-definable sets.)

We may now identify these sets of points with some familiar model theoretic 
objects. Let $\p$ be any partial type. The definable sets comprising it form 
a co-filtering system, with all maps the inclusions. The last equations says 
that computing the $M$ points of $\p$, viewed as pro-definable set, coincides 
with computing its $M$ points as a type, i.e., taking the intersection of the 
$M$ points of the definable sets in $\p$. The fact that two such system that 
give the same pro-definable set also give the same set of points means that 
this set of points is determined by the set of definable sets containing 
$\p$.

A partial type such as above is always contained in some definable set. There 
is a more general construction, called a $*$-type, that consists of the 
intersection of formulas in an arbitrary set of variables. Such types are 
similarly examples of pro-definable sets.

Analogously, an increasing union of definable set is an example of an 
ind-definable set. A more complicated example can be formulated as follows: 
let $E_i$ be definable equivalence relations on a definable set $X$, indexed 
by natural numbers $i$, such that for $i>j$, $E_i$ is coarser than $E_j$. Let 
$E$ be the equivalence relation saying that $xEy$ if $xE_i{}y$ for some $i$.  
Then $E$ is the union of the $E_i$ an thus an example of an ind-definable 
equivalence relation. The quotient of $X$ be $E$ is another example of an 
ind-definable set.

Our next purpose is to describe the morphisms between the new objects in 
terms of their points in models. Considering equations \eqref{eqn:homind} 
again, we see in particular that any morphism from $\Ind{(X_i)}$ to 
$\Ind{(Y_j)}$ gives rise to a filtering system $\Gamma_i$ of the 
corresponding graphs of functions from $X_i$ to $Y_{j_i}$. Similarly, a 
morphism of pro-definable sets gives rise to a cofiltering system. Each such 
system is isomorphic to its domain $X_i$, and therefore induces a function on 
the level of points from $\Ind{(X_i)}(M)$ to $\Ind{(Y_i)}(M)$ (and similarly 
for pro-definable sets.) We would like to show that conversely, any 
ind-definable set that gives rise to a function on the points of every model 
(equivalently, saturated enough model) induces a morphism.

We first restate the compactness theorem in this language:

\begin{prop}\label{prp:compact}
Let $\kappa$ be a cardinal bigger than the cardinality of the index category 
(i.e., the cardinality of the disjoint union of the morphism sets.)
\begin{enumerate}
  \item
    Let $f:\Ind{X_i}\ra{}Y$ be a morphism such that for some 
    $\kappa$-saturated model $M$, $f_M:\Ind{X_i}(M)\ra{}Y(M)$ is a bijection.  
    Then $f$ is an isomorphism.

  \item
    Let $f:Y\ra\Pro{X_i}$ be a morphism such that for some $\kappa$-saturated 
    model $M$, $f_M:Y(M)\ra\Pro{X_i}(M)$ is a bijection. Then $f$ is an 
    isomorphism.

\end{enumerate}
\end{prop}
\begin{proof}
  In each case, let $f_i$ be the maps corresponding to the morphism $f$. Note 
  that for definable sets and maps, the claims are true by definition ($f$ is  
  an isomorphism in this case.) We shall use the criterion of 
  remark~\ref{rmk:indiso} (and remark~\ref{rmk:proiso}.)

  \begin{enumerate}
    \item
      We will find $g$ and $t_i$ as required by remark~\ref{rmk:indiso}. We 
      first show that for some $k$, $f_k$ is onto. In fact, the collection of 
      sets $f_i(X_i)(M)$ is a small covering of $Y(M)$, hence it has a finite 
      sub-cover. Since the system is filtering, there is an $X_k$ above all 
      the sets in the sub-cover.
      
      We next note that the $t_i$ condition requires, in this case, for each 
      $i$, a definable map $t_i:X_i\ra{}X_j$ in the system such that 
      $f_i(x)=f_i(y)$ defines the same set as $t_i(x)=t_i(y)$. This again 
      holds by compactness: consider the set of formulas consisting of the 
      formula $f_i(x)=f_i(y)$, and for each morphism $t:X_i\ra{}X_j$ in 
      $(X_i)$, $t(x)\ne{}t(y)$. This set expresses the fact that the elements 
      $x,y\in{}X_i$ determine distinct elements of $\Ind{(X_i)}(M)$, that 
      have the same image under $f$. Therefore it is not satisfied in $M$. 
      Since this collection is small, a finite subset is not satisfied. 
      Therefore, there is some $t_i:X_i\ra{}X_j$ such that $f_i(x)=f_i(y)$ 
      implies $t_i(x)=t_i(y)$.

      In particular, this means that $f_j$ is injective on the image of 
      $t_i$. Let $X_0$ be the codomain of $t_k$ (for the $X_k$ found above.) 
      Then $f_0$ restricted to the image of $t_k$ is a bijection. We take $g$ 
      to be the inverse of this restriction.

    \item
      The proof is dual, using remark~\ref{rmk:proiso}. The only complication 
      here is that the category of definable sets does not, in general, have 
      finite direct limits. The assumption that such limits exist is called 
      \Def{elimination of imaginaries (EI)}. However, for the specific 
      purpose of the condition in remark~\ref{rmk:proiso}, we do not actually 
      need these limits. In our case, the condition simply translates to 
      saying that $f_i$ and $t_i$ have the same image. Such $t_i$ can be 
      obtained by compactness, using the surjectivity of the limit map, as in 
      the dual case.
      
      Further, by considering the formulas $f_i(x)\ne{}f_i(y)$, we see that 
      there is an $X_k$ such that $f_k$ is injective. In particular, we have 
      $t_k:X_0\ra{}X_k$, such that $f_k$ is a bijection between $Y$ and the 
      image of $t_k$. Taking $g=\inv{f_k}\circ{}t_k$, all the conditions of 
      the lemma are satisfied.
      
  \end{enumerate}
\end{proof}


The promised description of morphisms is just the extension of this criterion 
to the entire category:

\begin{prop}\label{prp:morphisms}
Let $\kappa$ be a cardinal bigger than the cardinality of the index category, 
$M$ a $\kappa$ saturated model. Let $X$ and $Y$ be ind- (or pro-) definable 
sets, $f:X\ra{}Y$ a morphism that induces a bijection on $M$ points. Then $f$ 
is an isomorphism.

In particular, there is a natural bijection between $\Hom{X}{Y}$ and 
sub-objects of $X\x{}Y$ whose set of $M$-points is a function from $X(M)$ to 
$Y(M)$.
\end{prop}
\begin{proof}
  We prove for the $Ind$ category, the $Pro$ case is dual. We have 
  $f:\Ind{X_i}\ra\Ind{Y_i}$. We first note that for any map $f:P\ra\Ind{Y_i}$ 
  where $P$ is a presheaf, $f$ is an isomorphism if and only if for all $j$, 
  the pullback $f_j:P\x_{\Ind{Y_i}}Y_j\ra{}Y_j$ is an isomorphism. Indeed, 
  given inverses $g_j$ to the $f_j$, their composition with the projection to 
  $P$ forms a matching family of maps from the $Y_i$ to $P$, and therefore 
  yields a map from $\Ind{Y_i}$ to $P$, inverse to $f$.

  Furthermore, if $P$ itself is ind-definable, $P=\Ind{X_i}$, we have
  \begin{equation*}
    P\x_{\Ind{Y_i}}Y_j=\Ind{X_i\x_{\Ind{Y_i}}Y_j}=\Ind{X_i\x_{Y_{k_i}}Y_j}
  \end{equation*}
  On the other hand, since taking $M$ points is represented by a 
  pro-definable set, it preserves pullbacks. Therefore, if $f_M$ is a 
  bijection of $M$ points, so is ${f_j}_M$, for any $j$. By 
  proposition~\ref{prp:compact}, $f_j$ is an isomorphism.

  The description of the morphism sets is the interpretation of this 
  statement for the projection map from a sub object $R$ of $X\x{}Y$ to $X$.

\end{proof}

We may summarise the results of this section as follows:

\begin{cor} \label{cor:proind}
  Let $M$ be a $\kappa$-saturated model.

  The functor of ``taking $M$ points'' is an equivalence of categories 
  between the category $\Pro[\kappa]{\DD}$ of pro-definable sets 
  representable by systems of length less than $\kappa$, and the sub-category 
  of the category of sets whose objects and morphisms are inverse co-filtered 
  limits of $M$ points of definable sets, of length less than $\kappa$.

  Similarly, the same functor is an equivalence of categories between the 
  category $\Ind[\kappa]{\DD}$ of ind-definable sets representable by systems 
  of length less than $\kappa$, and the sub-category of the category of sets 
  whose objects and morphisms are direct filtered limits of $M$ points of 
  definable sets, of length less than $\kappa$.

\end{cor}

Finally, we note that definable sets are given with canonical inclusions (in 
the ``universe''.) For example, in our terminology, any two points are 
identified. If we wish to remember the inclusion of the definable sets in 
some definable set $X$, we work in the category $\DD/X$, and all results 
continue to hold. This way we get pro-definable subsets of $X$. These sets 
are called also \Def{$\w$-definable}.

\bibliography{../../bibtex/mr}

\end{document}